\documentclass[a4paper,10pt,reqno,twoside]{amsart}

\usepackage{mathtools,amsmath,amssymb,amsfonts,amsthm,mathrsfs}
\usepackage[pdftex]{graphicx}
\usepackage{graphicx}
\mathtoolsset{showonlyrefs}

\DeclareGraphicsExtensions{.bmp,.png,.pdf,.jpg,}

\usepackage{bm}
\usepackage{verbatim}
\usepackage{color}
\usepackage{mathrsfs}
\usepackage{pgfplots}
\usepackage{euscript}
\usepackage{bbm}
\usepackage{cite}

\newcommand{\RR}{\mathbb R}

\newcommand{\sgn}{{ \rm sgn }}

\newcommand{\limess}[1]%
{

\begin{array}[t]{c}
{\rm ess\, lim}\\
{\scriptstyle #1}
\end{array}

}

\newcommand{\supess}[1]%
{

\begin{array}[t]{c}
{\rm ess\, sup}\\
{\scriptstyle #1}
\end{array}

}
\newcommand{\infess}[1]%
{

\begin{array}[t]{c}
{\rm ess\, inf}\\
{\scriptstyle #1}
\end{array}

}

\newcommand{\II}{{\mathcal I}}

\def\R{{\mathbb R}}

\def\a{\alpha}

\def\d{\delta}

\def\la{\lambda}

\def\O{\Omega}
\def\e{\varepsilon}

\def\o{\omega}

\def\mf{\mathfrak}

\numberwithin{equation}{section}

\theoremstyle{definition}
\newtheorem{dfn}{Definition}[section]

\theoremstyle{definition}
\newtheorem{remark}{Remark}[section]
\newtheorem{example}{Example}[section]

\theoremstyle{plain}
\newtheorem{theorem}{Theorem}[section]

\newtheorem{lemma}{Lemma}[section]

\pgfplotsset{compat=1.17}

\begin{document}

\title[Regularity results    via the strong maximum principle]
{Optimal regularity  results for   the one-dimensional  prescribed curvature equation via the strong
maximum principle
\vskip1cm}

\thanks{J. L\'{o}pez-G\'{o}mez has been  supported by the Research Grants PGC2018-097104-B-I00 and 
PID2021-123343NB-I00 of the Spanish Ministry of Science and Innovation, and by the Institute of Interdisciplinary Mathematics of Complutense University.
\\
\indent
This research has been performed under the auspices of INdAM-GNAMPA}

\author{Juli\'{a}n L\'{o}pez-G\'{o}mez}
\address{Juli\'{a}n L\'{o}pez-G\'{o}mez:
Instituto Interdisciplinar de Matem\'aticas,
Universidad Complutense de Madrid,
Madrid, Spain}
\email{jlopezgo@ucm.es}
\author{Pierpaolo Omari}
\address{Pierpaolo Omari:
Dipartimento di Matematica e Geoscienze,
Universit\`a degli Studi di Trieste,
Trieste, Italy}
\email{omari@units.it}

\begin{abstract}
\vskip.5cm
A refined version  of the strong maximum principle  is   proven for a  class of second order ordinary differential equations
with  possibly discontinuous non-monotone nonlinearities.
Then, exploiting this tool,   some optimal regularity results  recently established by L\'opez-G\'omez and Omari,
in \cite{LGO-JDE20, LGO-ANS20, LGO-21},
for the   bounded variation  solutions of non-autonomous quasilinear equations driven by the  one-dimensional curvature operator, are
 substantially improved
by admitting general prescribed curvatures and by incorporating  general boundary conditions.
The new approach developed here yields a new, deeper,    interpretation of the  assumptions  introduced in our  previous papers, simultaneously clarifying their meaning and making fully transparent their connection  with the strong maximum principle.
\vskip.5cm
\end{abstract}


\vspace{0.1cm}

\date{\today}

\maketitle

\section{Introduction}
\label{s1}

\noindent The aim of this work is twofold. First we prove an extended version  of the   strong maximum principle for  a general class of  second order ordinary differential equations
\begin{equation}
\label{g}
v'' =g(t, v, v'),
\end{equation}
in the absence of any assumption of continuity or monotonicity,
which is a crucial feature at the light of our subsequent applications.
Then, based on this  form of the strong maximum principle,  we provide   some
optimal regularity results for the bounded variation solutions, positive and nodal,  of the non-autonomous    curvature  equation
\begin{equation}
\label{e}
-\left( \frac{u'}{\sqrt{1+(u')^2}}\right)' = f(x,u),
\end{equation}
where $f$ is an arbitrary function prescribing the curvature of $u$. These  findings substantially generalize some previous regularity results we established in \cite{LGO-ANS19, LGO-JDE20, LGO-ANS20, LGO-21} for the positive bounded variation solutions of \eqref{e}, under homogeneous Neumann boundary conditions, in the special case where
\begin{equation}
\label{fhk}
f(x,s) = h(x)\, k(s).
\end{equation}
The reader should be aware that, unlike in these papers, here we are   imposing neither the structural assumption \eqref{fhk}, nor the Neumann boundary conditions, nor any requirement on the sign  on the solutions.  As a consequence, the analysis  carried out here allows us, through a completely different technical device,  to extend  most of the results of   \cite{LGO-ANS19, LGO-JDE20, LGO-ANS20, LGO-21}  to  more  general classes of equations and boundary value problems. Furthermore, in this paper we are going to provide with a new   interpretation of  the assumptions  used in  our previous works,  clarifying their meaning and establishing some deep,   though previously hidden, connections with the strong maximum principle.

This paper is organized as follows.
In Section \ref{s2} we establish  as Theorem \ref{th2.1} the version of the strong maximum principle for equation \eqref{g}  used  throughout this paper. Since
no  kind  of continuity or monotonicity,  even  in the state variable $s$,  is imposed
on the right hand side $g= g(t,s,\xi)$ of equation \eqref{g}, as well as  on the comparison function $G'= G'(s)$ introduced in condition (G) of
Theorem \ref{th2.1}, this form of the   strong maximum principle  provides a completion and a sharpening  of  its  counterparts in \cite{Va84} or \cite{PuSe}; its proof being also  more
delicate than in the classical situations. 

Section  \ref{s3}  contains the main regularity results of this work. The first one, stated as Theorem \ref{th3.1}, shows  that a bounded variation solution $u$ of the curvature equation \eqref{e}
can loose its regularity  only at the endpoints, but never
at the interior points,
of  any interval where  the function $f(\cdot, u(\cdot))$  has a definite sign; yet,  $u$ can be singular at an interior point of its domain if such a point separates two adjacent intervals where $f(\cdot, u(\cdot))$ changes sign. This result sharpens \cite[Prop. 3.6]{LGO-ANS19}, where the  solutions were assumed, in addition, to satisfy homogeneous  Neumann boundary conditions.
Although our  proof of  Theorem \ref{th3.1} is a re-elaboration of the proof of \cite[Prop. 3.6]{LGO-ANS19}, for the sake of completeness, we are giving complete technical details here.

Our next two results, Theorems \ref{th3.2}  and \ref{th3.3},  establish the complete  regularity of the bounded variation solutions $u$ of  \eqref{e}. Precisely,  Theorem \ref{th3.2}  guarantees the  regularity at the endpoints of any interval where  the sign of $f(\cdot, u(\cdot))$ is constant, by imposing at these points a suitable control, expressed by any of the conditions (j)--(jjjj),  on  the  decay rate to zero of $f(\cdot, u(\cdot))$
Theorem \ref{th3.3},  instead, guarantees the regularity of $u$ at any interior point, $z$, separating two adjacent interval where $f(\cdot, u(\cdot))$ changes sign, by imposing a similar decay property to $f(\cdot, u(\cdot))$ either on  the left, or on the right, of $z$, as expressed   by the conditions (h) or (hh).

 Such  kind of decay controls  were introduced by the authors in \cite{LGO-JDE20} for discussing the regularity  of the positive   bounded variation   solutions of the Neumann problem associated with   \eqref{e} in the  very special case when  $f$ can be decomposed as in \eqref{fhk}, namely,
 $ f(x,s) = h(x)\, k(s), $
 where $k$ is a function having a superlinear potential at infinity. They were later used by the authors in   \cite{LGO-ANS20}   to treat the case where the  potential of $k$ is asymptotically linear, and  in \cite{LGO-21} to deal with the case of sublinear potentials. Here, no specific restriction on the asymptotic behavior of $f(\cdot,s)$  with respect to $s$ is  imposed. From \cite{LGO-JDE20, LGO-ANS20, LGO-21} we also  know  that these  assumptions on the decay rate of $f(\cdot, u(\cdot))$ are optimal, in the sense that, if  they fail at some point, the derivative $u'$ might blow-up  there, and the solution $u$ might even develop a jump discontinuity.

The proof of Theorems \ref{th3.2} and \ref{th3.3} presented here is completely new and it relies on the use of the strong maximum principle as expressed by Theorem \ref{th2.1}.
Our   approach, besides being far more general and versatile, displays the following striking fact. It turns out that the precise  condition yielding the regularity  of a solution $u$ of  \eqref{e}, through a control on the decay rate to zero of $f(\cdot, u(\cdot))$ at some point $z$,
is precisely  the assumption required by Theorem \ref{th2.1} so that the strong maximum principle holds for  the differential  equation
\begin{equation}
\label{inv}
\left( \frac{v'}{\sqrt{1+(v')^2}}\right)' =
 f(z+v,t)
\quad \Longleftrightarrow \quad
v''= f(z+v,t)\big(1+(v')^2\big)^\frac{3}{2},
\end{equation}
satisfied by the shift $v= w-z$ of a  local inverse $w$ of $u$. Note that, as $f$
is not assumed to satisfy
any regularity condition, the right hand side of \eqref{inv}, that is, the function
$$
g(t,s,\xi) := f(z+ s, t) \big(1+\xi^2\big)^\frac{3}{2},
$$
may be discontinuous, besides in $t$, in the state variable $s$ as well. Note that this could happen for $g$ even if $f$ were a Carath\'eodory function.

Essentially, we establish that the validity of the strong maximum principle for equation \eqref{inv} yields the regularity for the solutions of \eqref{e}. As a consequence, the bounded variation solutions of \eqref{e} can develop singularities only when the conclusions of the strong maximum principle fail for \eqref{inv}.
This appears to be a quite remarkable achievement that
illuminates and clarify the apparently  exotic conditions introduced in \cite{LGO-JDE20}.

Finally, Section \ref{s4} is devoted to the application of the  theory developed in Section \ref{s3} to
establish the regularity, up to the boundary, of the solutions of a number of non-autonomous one-dimensional prescribed curvature equations supplemented with several types of boundary conditions, such as Dirichlet, Neumann,   Robin, or even periodic boundary conditions. These and other similar statements, that can be deduced from Theorems \ref{th3.1}, \ref{th3.2},  \ref{th3.3},  complete or extend, as far as  regularity is concerned, several  existence and multiplicity results  previously obtained   in \cite{BHOO-JDE, BHOO-TS, BOO-DIE, CLT-JDE, CLT-NA, OO-DIE, OO-JDE10, OO-CM, OO-DCDS, OS-NA, LOR1, LOR2, LGO-ANS19, LGO-JDE20, LGO-ANS20, LGO-21, LGO-AML22, OO-JDE11, OO-CCM, OOR-NARWA, OOR-JFA, COZ-CCM, OO-OM}.

\section{A variant of the strong maximum principle}
\label{s2}

\noindent The main result of this section   is  the following       version of the strong maximum principle for second order ordinary differential equations  with  possibly discontinuous non-monotone  right hand sides.

\begin{theorem}
\label{th2.1}
Let $g:  (\a, \o) \times \RR \times \RR \to \RR$  be a given function  and let $v\in W^{2,1}_{\rm loc}(\a,\o) \cap W^{1,1}(\a,\o)$  be a non-trivial non-negative solution of the  differential equation
\begin{equation}
\label{v''g}
  v''(t)  =g(t, v(t), v'(t)) \quad  \hbox{for almost all}\;\; t\in (\a, \o).
\end{equation}
Assume that:
\smallskip
\begin{itemize}
\item[{\rm (G)}]
there exist a constant $\e>0$ and an absolutely continuous function $G: [0,\e] \to \RR$ such that
\begin{equation}
\label{0gG'}
0\le g(t,v(t), v'(t)) \le G'(v(t)) \;\; \hbox{for almost all}\;\; t\in (\a, \o)\;\;\hbox{for which}\;\;
 0<v(t) \le \e\;\;\hbox{and}\;\; |v'(t)|\le \e,
\end{equation}
and either
\begin{equation}
\label{G=0}
G(s) =0 \;\; \hbox{for all} \;\; s \in (0,\e],
\end{equation}
or
\begin{equation}
\label{G>0}
   G(s)>0 \;\; \hbox{for all} \;\; s \in (0,\e] \;\;\hbox{ and }\;\; \int_0^\e \frac{1}{\sqrt{G(s)}} ds = +\infty.
\end{equation}
\end{itemize}
Then,  $v$ is strongly positive, in the sense that  the following properties hold true:
\begin{itemize}
\item[{\rm (i)}] $v(t) >0$ for  all $t\in (\a, \o)$,
\smallskip
\item[{\rm (ii)}]  $v'(\a^+)>0$ \hbox{if} $v(\a) =0$ and $v'(\a^+)$ exists,
\smallskip
\item[{\rm (iii)}] $v'(\o^-)<0$ if $v(\o) =0$ and  $v'(\o^-)$   exists.
\end{itemize}
\end{theorem}

\begin{remark}
It is worth observing that no kind of continuity, or monotonicity, is  imposed either on  the function $g$ appearing at the right hand side of equation \eqref{v''g}, or on the comparison function $G'$.
Indeed,  the function $G'$ considered in  assumption (G)
is   only Lebesgue integrable and does not  satisfy any monotonicity condition. These features   are crucial  for proving our new findings in Section \ref{s3}.
\par
In this respect, assumption (G) is independent of the  conditions required by the classical V\'azquez strong maximum principle in \cite{Va84}, where $G'$ is assumed to be continuous and increasing,  as well as of  the conditions  imposed  by Pucci and Serrin in \cite[Ch. 5]{PuSe}, where $G'$ is again  supposed  to be increasing. Th
\par
As already  pointed out in \cite{Va84,  PuSe},   the condition \eqref{G>0} is sharp, because if it fails, dead core solutions, i.e., non-negative solutions vanishing on sets of positive measure, may occur.    The necessity of this type of conditions goes back to Benilan, Brezis and Crandall \cite{BBC}.
Note that condition (G) entails $G(0)=0$.

\end{remark}

\begin{remark}
Under condition (G),   there are some positive constants, $k$, that are strict supersolutions of the
differential operator
$$
  \mathfrak{L}v:= -v''+G'(v), \qquad v\in W^{2,1}_{\rm loc}(\a,\o) \cap W^{1,1}(\a,\o),
$$
under Dirichlet boundary conditions in $(\a,\o)$.
Thus, $\mathfrak{L}$ satisfies similar assumptions as those of Theorem 10 in Chapter 2 of Protter and Weinberger \cite{PW}. From such a perspective, Theorem \ref{th2.1} can be viewed as a
sharp
 one-dimensional nonlinear counterpart of Corollary 2.1 of L\'{o}pez-G\'{o}mez \cite{LG-WS} (going back to \cite{LG-ADE} and \cite{LG-RIMS}). Indeed, one might re-state Theorem \ref{th2.1} by simply saying that, under condition (G), any  non-trivial non-negative solution 
of $\mathfrak{L}v\geq 0$ in $W^{2,1}_{\rm loc}(\a,\o) \cap W^{1,1}(\a,\o)$ must be strongly positive.
\par
This bi-association suggests the validity of Theorem \ref{th2.1}, even in a multidimensional context, for general operators of the form $\mathfrak{L}+G'(v)$ where $\mf{L}$ is a linear second order uniformly elliptic operator in $\O$ whose principal eigenvalue, under Dirichlet boundary conditions in $\O$, is positive.
\end{remark}

The proof of Theorem \ref{th2.1} exploits  the following maximum principle  for   first order   differential inequalities.

 \begin{lemma}
\label{le2.1}
Assume that $H:\RR \to \RR$
is a Lebesgue measurable function satisfying the condition:
\begin{itemize}
\item[{\rm (H)}]  there exists a constant $\e>0$ such that $H(s)> 0$ for all $s~\in [-\e, \e] \setminus \{0\}$ and
\begin{equation}
\label{infty}
\int_{-\e}^0 \frac{1}{H(s)}\,ds  = + \infty= \int_0^\e \frac{1}{H(s)}\,ds.
\end{equation}
\end{itemize}
Then, any non-trivial solution $v\in W^{1,1}(\a, \o) $
of the differential inequality
\begin{equation}
\label{v'H}
   |v'(t) | \le H(v(t)) \quad \text{for almost all } t\in (\a, \o)
\end{equation}
 is either strictly positive,  or strictly negative, i.e.,  either $\min v >0$,   or $\max v<0$.
\end{lemma}

\begin{remark}
Assumption (H) is the classical Osgood condition, introduced  in \cite{Osg} to guarantee  the uniqueness of the solution for the Cauchy problem associated with first order ordinary differential equations.
\end{remark}

The proof of Lemma  \ref{le2.1}   is elementary when $H$ is continuous. As here we are
only imposing  measurability, our proof is far more delicate.

\begin{proof}[\bf Proof of Lemma \ref{le2.1}]
Let $v\in W^{1,1}(\a, \o) $ be a non-trivial solution of \eqref{v'H}.   We claim that
$\min v >0$ if $\max v >0$.  It is apparent that,   substituting $-v$  for  $v$,
we can infer that  $\max v <0$ if $\min v <0$.
Arguing by contradiction,  we suppose    that
$$
  \max v >0 \quad \text{and} \quad \min v \le 0.
$$
As $v$ is continuous and $\max v >0$, there exist a point  $t_0\in [\a, \o] $ and two constants $\d>0$ and $\eta \in (0,\e]$ such that $v(t_0)=0$ and either
$$
v((t_0, t_0+\d])=(0, \eta],  \quad \text{or} \quad
v([t_0-\d, t_0))=(0, \eta].
$$
Assume that the former case occurs, the latter one being teated similarly.
As $H$ is   measurable   and $H(s)>0$ for all $s\in (0,\eta]$, the function  $\frac{1}{H}: (0, \eta] \to \RR$  is well defined and    it is   measurable too.
Moreover,   being absolutely continuous, $v$  satisfies the Lusin's $N$-property, i.e., it maps sets of null    measure  to sets of null   measure. Thus, the functions
$$
\frac{1}{H}: (0, \eta] \to \RR\quad \text{and} \quad  v:(t_0, t_0+\d]\to  (0, \eta]
$$ fulfill  the assumptions of \cite[Thm.2]{Ha}. Consequently,  the functions
$$
  \frac{|v'|}{H\circ v} : (t_0, t_0+\d]\to\RR\quad \hbox{and}\quad \frac{N_v(\cdot,(t_0,t_0+\d])}{H(\cdot)}  : (0, \eta] \to \RR
$$
are    measurable,   where
$
N_v(s,(t_0,t_0+\d])
$
denotes the Banach indicatrix  of $v$ in the interval $(t_0,t_0+\d]$, i.e.,
$$
N_v(s,(t_0,t_0+\d]) = \mathcal{H}^0\left( v^{-1} (s) \cap (t_0,t_0+\d] \right),
$$
with  $\mathcal{H}^0$   the counting measure.
Moreover, by \eqref{v'H}, we have that
\begin{equation}
\frac{|v'(t)| }{ H(v(t)) } \le 1 \quad  \text{for almost all } t\in (t_0,t_0+\d].
\end{equation}
Hence,    integrating yields
\begin{equation}
\label{led}
 \int_{t_0}^{t_0+\d}\frac{|v'(t)| }{ H(v(t) }\, dt \le \d.
\end{equation}
 Furthermore, by the  formula of change of variables established in \cite[Thm.2]{Ha}, we find that
\begin{equation}
\label{cv}
 \int_{t_0}^{t_0+\d}\frac{|v'(t)| }{ H(v(t)) }\, dt = \int_0^\eta \frac{N_v(s,(x_0,x_0+\d])}{H(s)}  \, ds.
 \end{equation}
As $v((t_0, t_0+\d])=(0, \eta]$,  we have that
$$
N_v(s,(t_0,t_0+\d]) \ge 1 \quad \text{for all } s \in (0, \eta].
$$
 Hence,  by \eqref{cv} and \eqref{led}, we get
 \begin{equation}
 \int_0^\eta \frac{1}{H(s)}  \, ds \le \int_0^\eta \frac{N_v(s,(x_0,x_0+\d])}{H(s)}  \, ds= \int_{t_0}^{t_0+\d}\frac{|v'(t)| }{ H(v(t)) }\, dt \le \d,
  \end{equation}
which is impossible, because \eqref{infty}  entails  that
$$
 \int_0^\eta \frac{1}{H(s)}  \, ds = +\infty\quad \text{for all } \eta \in (0,\e].
$$
This concludes the proof of Lemma \ref{le2.1}.
\end{proof}

Similarly as the proof of Lemma \ref{le2.1} also the proof of Theorem \ref{th2.1} would be      a bit simpler if, in addition, $g$ and $G'$ are continuous and   $v\in C^2[\a,\o]$.

\begin{proof}[\bf Proof of Theorem \ref{th2.1}]
Let $v\in W^{2,1}_{\rm loc}(\a,\o) \cap W^{1,1}(\a,\o)$  be a non-trivial non-negative solution of  \eqref{v''g}.
Arguing by contradiction, we suppose that there is a point $t_0\in [\a,\o]$ such that
\begin{equation}
\label{v=0}
 \min v = v(t_0) =0.
\end{equation}
Obviously,  we must have $v'(t_0)=0$, if $t_0\in (\a,\o)$.
 We further  assume that  $v'(t_0)= v'(\a^+)=0$, if $t_0=\a$, or respectively $v'(t_0)= v'(\o^-)=0$, if $t_0=\o$, provided  that $v'(\a^+) $, or  $v'(\o^-) $, exists. As $\max v>0$, we can also suppose, possibly for a different choice of $t_0$, that there exists a constant $\d>0$ such that
$$
0<v(t) \le \e \quad \text{and} \quad |v'(t)| \le \e,
$$
either for all $t\in (t_0,t_0+\d]$, or for all $t\in [t_0+\d,t_0)$, where $\e$ is the constant appearing in condition (G).
Assume that the former case occurs.  The proof, being similar in the latter one, will be omitted here. From \eqref{v''g} and \eqref{0gG'} it follows  that $v''(t) \ge 0$   for almost all  $t\in (t_0,t_0+\d]$. Moreover, since
$$
v'(t) - v'(s)= \int_{s}^t v''(r) \, dr \quad \text{for every $s, t\in (t_0,t_0+\d]$,}
$$
 letting $s\to t_0$, we infer that
$$
v'(t)  = \int_{t_0}^t v''(r) \, dr \quad \text{for every $t\in (t_0,t_0+\d]$.}
$$
Thus,  $v' $ is absolutely continuous and increasing in $[t_0,t_0+\d]$.
Moreover, as $v(t) >0$ for all $t\in (t_0,t_0+\d]$ and
 $v(t_0)=0=v'(t_0)$, we conclude that
 \begin{equation}
 \label{u'>0}
v'(t) >0 \quad \text{for every } t\in (t_0,t_0+\d].
\end{equation}
From \eqref{0gG'},  we have that
$$
   g(s,v(s), v'(s)) \le G'(v(s)) \quad \hbox{for almost all}\;\; s\in [t_0,t_0+\d],
$$
and hence, by \eqref{v''g},
\begin{equation}
\label{u''G}
v''(s) \le G'(v(s)) \quad   \hbox{for almost all}\;\; s\in [t_0,t_0+\d].
\end{equation}
Consequently, multiplying \eqref{u''G} by $v'(s)$, it follows from \eqref{u'>0}  and \eqref{u''G} that
\begin{equation}
\label{u''u'Gu'}
v'(s)\,  v''(s)  \le  G'(v(s)) \, v'(s) \quad \hbox{for   almost all}\;\; s\in [t_0,t_0+\d].
\end{equation}
Since $v$ is absolutely continuous   and  strictly increasing in $ [t_0,t_0+\d]$ and $G$ is absolutely continuous in $[0,v(t_0+\d)]$,    the composition $G\circ v$ is absolutely continuous    in $ [t_0,t_0+\d]$. Clearly,  $(v')^2 $ is absolutely continuous    in $ [t_0,t_0+\d]$ too. Therefore, for every $t\in (t_0,t_0+\d]$, integrating \eqref{u''u'Gu'} in $[t_0, t]$ and  applying the formula of change of variables yields
\begin{equation}
\label{u'G}
 \frac{1}{2}(v'(t))^2  = \int_{t_0}^{t } v'(s) v''(s) \, ds \le \int_{t_0}^t G'(v(s)) v'(s) \, ds =  \int_{v(t_0) }^{v(t)} G'(s)  \, ds  = G(v(t)),
\end{equation}
because  $v'(t_0)=0$ and $G(v(t_0))=G(0) =0$.
 From \eqref{u'G} it follows that
 \begin{equation}
|v'(t)| \le   \sqrt{2G(v(t))} \quad\hbox{for every}\;\; t \in   [t_0,t_0+\d].
\end{equation}
In case  \eqref{G=0} holds, we find that
$$
v'(t)=0 \quad \text{for all } t \in   (t_0,t_0+\d],
$$
thus contradicting \eqref{u'>0}.   Whereas, if \eqref{G>0} holds, then the function  $H= \sqrt{2G}$ satisfies the assumptions  of Lemma \ref{le2.1}, which  implies that
$$
 v(t) >0 \quad \text{for all } t \in   [t_0,t_0+\d],
 $$
  thus contradictig \eqref{v=0}. This concludes the proof of Theorem \ref{th2.1}.
\end{proof}

\section{Optimal regularity results for the  prescribed curvature equation}
\label{s3}

\noindent In this section we discuss the regularity properties of the bounded variation solutions of the one-dimensional non-autonomous prescribed  curvature equation
 \begin{equation}
 \label{E}
-\left( \frac{u'}{\sqrt{1+(u')^2}}\right)' = f(x,u), \qquad a<x<b,
\end{equation}
where  $f:(a,b)\times\RR \to \RR$  is any  given  function.
  We begin by recalling   the notion  of  bounded variation solution of equation \eqref{E}.
To this end, for any $v\in BV(a,b)$, we denote by $ Dv =  {D^a v}\,  dx +  {D^s v}$  the Lebesgue--Nikodym decomposition, with respect to the Lebesgue measure $dx$ in $\RR$,  of the Radon measure $Dv$ in its absolutely continuous part $  {D^av}\, dx$, with density function ${D^av}$,  and  its singular part  ${D^s v} $. Further,  $ \frac{ {D^sv}}{ |{D^sv}|}$  stands for the density function of  ${D^sv}$ with respect to its absolute variation $|{D^sv}|$.  Finally,   for every  $x_0\in [a, b)$, $v(x_0^+) $ denotes the  right trace of $v$ at $x_0$ and,   for every  $x_0\in (a, b]$,  $v(x_0^-) $ denotes the  left trace of $v$ at $x_0$.  We refer, e.g.,  to  \cite{AmFuPa} for additional details on these concepts.

\begin{dfn}
A function  $u\in BV(a,b)$ is a bounded  variation solution of  \eqref{E}
if $f(\cdot,u(\cdot)) \in L^1(a,b)$  and
\begin{equation}
\label{EE}
\int_a^b \frac{{D^au}(x) D^a\phi(x) }{\sqrt{1+(D^au(x))^2 }} \, dx + \int_a^b  \frac{D^su}{|D^su|} (x) \,    D^s\phi = \int_a^b  f(x,u(x)) \phi(x) \, dx
\end{equation}
for all $\phi \in BV(a,b) $ such that $|D^s\phi| $ is absolutely continuous with respect to $ |D^su|$ and $\phi(a^+)= \phi(b^-)=0$.
\end{dfn}

\begin{remark}
It is proven in   \cite{Anz83} that $u\in BV(a,b)$ is a bounded  variation solution of  \eqref{E}  if and only if it minimizes in $BV(a,b)$ the functional
\begin{equation}
\label{M}
\II(v) := \int_a^b \sqrt{1+| D v|^2}  + |v(a^+) -u(a^+)| + |v(b^-) -u(b^-)| - \int_a^b f(x,u(x)) v(x) \, dx,
\end{equation}
where
$$
\int_a^b \sqrt{1+| D v|^2} = \int_a^b \sqrt{1+| D^av(x)|^2}\, dx  +  \int_a^b   |D^sv|.
$$
\end{remark}

The next regularity result establishes that a bounded variation solution $u$ of \eqref{E} can loose its regularity   at the endpoints, but never at the interior points,  of  the intervals where  the function $f(\cdot,u(\cdot) )$  has a definite sign; whereas, $u$  can be singular at an interior point of its domain if such a point separates two adjacent intervals where $f(\cdot,u(\cdot) )$ changes sign. In both cases, the derivative $u'$ blows up, but, in the latter one, $u$ can further exhibit a jump discontinuity.
Although Theorem \ref{th3.1} is a substantial refinement of \cite[Prop. 3.6]{LGO-ANS19}, which was limited to the  solutions of   equation \eqref{E} satisfying the Neumann boundary conditions $u'(a)= u'(b)=0$, its proof is basically a re-elaboration of the proof of \cite[Prop. 3.6]{LGO-ANS19}. Nevertheless, for the sake of completeness, we are delivering complete technical details here.

\begin{theorem}
\label{th3.1}
Let $u$ be a bounded variation solution of equation \eqref{E}. Then, the following statements hold.
\begin{itemize}
\item[${\rm (i})$]
If $ f(x,u(x))\ge 0 $   for almost all $x\in (a,b)$,  then $u$ is concave and
either $ u\in W^{2,1}(a,b)$,
or $ u\in W^{2,1}_{\text{\rm loc}}[a,b) \cap W^{1,1}(a,b) $ and $u'(b^-)=-\infty$,
or $ u\in W^{2,1}_{\text{\rm loc}}(a,b] \cap W^{1,1}(a,b)$ and $u'(a^+)=+\infty$,
or $ u\in W^{2,1}_{\text{\rm loc}}(a,b) \cap W^{1,1}(a,b)$, $u'(a^+)=+\infty$, and $u'(b^-)=-\infty$.
In all cases, $u$ satisfies   equation \eqref{E} for almost all $x\in (a, b)$.
\smallskip
\item[{\rm (ii)}]
If $ f(x,u(x))\le 0 $   for almost all $x\in (a,b)$,  then $u$ is convex and
either $ u\in W^{2,1}(a,b)$,
or $ u\in W^{2,1}_{\text{\rm loc}}[a,b)\cap W^{1,1}(a,b)$ and $u'(b^-)=+\infty$,
or $ u\in W^{2,1}_{\text{\rm loc}}(a,b]\cap W^{1,1}(a,b)$ and $u'(a^+)=-\infty$,
or $ u\in W^{2,1}_{\text{\rm loc}}(a,b) \cap W^{1,1}(a,b)$, $u'(a^+)=-\infty$, and $u'(b^-)=+\infty$.
In all cases, $u$ satisfies   equation \eqref{E} for almost all $x\in (a, b)$.
\smallskip
\item[{\rm (iii)}]
If there is $c\in (a,b)$ such that $ f(x,u(x))\ge 0 $   for almost all $x\in (a,c)$ and $ f(x,u(x))\le 0 $   for almost all $x\in (c, b)$, then $u_{|(a,c)}$ is concave, $u_{|(c,b)}$ is convex, and
either $ u\in W^{2,1}_{\text{\rm loc}}(a,b) \cap W^{1,1}(a,b)$,
or $ u_{|(a,c)}\in W^{2,1}_{\text{\rm loc}}(a,c) \cap W^{1,1}(a,c)$, $u_{|(c,b)} \in  W^{2,1}_{\text{\rm loc}}(c,b)\cap W^{1,1}(c,b)$, $u(c^-) \ge  u(c^+)$, and $u'(c^-)=-\infty= u'(c^+)$. Moreover, in case $u(c^-) > u(c^+)$, we have that
\begin{equation}
\label{delta}
D^su = (u(c^+) - u(c^-))\, \delta_c,
\end{equation}
where $\delta_c$ stands for   the Dirac measure concentrated at $c$. In any circumstances,  $u$ satisfies equation \eqref{E} for almost all $x\in (a, b)$.
\smallskip
\item[{\rm (iiii)}]
If there is $c\in (a,b)$ such that $ f(x,u(x))\le 0 $   for almost all $x\in (a,c)$ and $ f(x,u(x))\ge 0 $   for almost all $x\in (c, b)$, then $u_{|(a,c)}$ is convex, $u_{|(c,b)}$ is concave, and
either $ u\in W^{2,1}_{\text{\rm loc}}(a,b) \cap W^{1,1}(a,b)$,
or $ u_{|(a,c)}\in W^{2,1}_{\text{\rm loc}}(a,c) \cap W^{1,1}(a,c)$, $u_{|(c,b)} \in  W^{2,1}_{\text{\rm loc}}(c,b)\cap W^{1,1}(c,b)$, $u(c^-) \le  u(c^+)$, and $u'(c^-)=+\infty= u'(c^+)$. Moreover, in case $u(c^-) < u(c^+)$, \eqref{delta} holds. In any circumstances, $u$ satisfies equation \eqref{E} for almost all $x\in (a, b)$.
\end{itemize}
\end{theorem}

\begin{remark}
In cases (iii) and (iiii), the behavior of $u$ at the endpoints
$a$ and   $b$
follows exactly the  same patterns as  described in (i) and (ii).
\end{remark}

\begin{proof}
Let $u$ be a bounded variation solution of   equation \eqref{E}, and consider the  decomposition of the measure $Du$,
$$
Du = D^au\,dx+ D^su = D^au\,dx+ D^ju + D^cu,
$$
in   its absolutely continuous part, $D^au\,dx$, its jump part, $D^ju$, and its Cantor part, $D^cu$, as well as  the induced decomposition   of the function $u$,
$$
  u = u^a + u^j + u^c,
$$
where
$$
   Du^a\,dx = D^a u\, dx,\quad Du^j = D^j u,\quad Du^c = D^c u,
$$
and
$$
  u^a(a) = u(a^+), \quad u^j(a^+) = 0, \quad u^c(a) = 0.
$$
Hereafter,  for convenience,  we write
$$
  v:=u^a\in W^{1,1}(a, b),\quad h := f(\cdot, u(\cdot))\in L^1(a,b), \quad \psi(s) := \frac{s}{\sqrt{1+s^2 }},\;\text{for all }s\in\RR.
$$
\smallskip
To prove the theorem, it suffices to show that the Assertions (i) and (iii) hold true.
\medskip

\noindent \textit{Proof of Assertion \rm (i).}
The proof is divided into three steps.
\smallskip

\noindent\texttt{Step 1}: ${u^a} \in W^{2,1}_{\text{{\rm loc}}}(a,b)$, and it is concave in  $(a, b)$.
\\
Testing \eqref{EE} against functions $\phi \in W^{1,1}(a,b)$ with $\phi(a)= \phi(b)=0$,  yields
$$
 \int_a^b
\psi (v'(x)) \phi'(x)\, dx   = \int_a^b  h(x) \phi(x) \, dx.
$$
 Thus, $\psi(v')\in  W^{1,1}(a, b)$ and
\begin{equation}
\label{eqpsi}
-\big(\psi(v'(x))\big)'= h(x) \quad \hbox{for almost all } x \in  (a, b).
\end{equation}
As $h(x) \ge 0 $  for almost all  $x \in (a, b)$,
$\psi(v') $ is decreasing in $(a, b)$.  Since, in addition, $\psi(v') $ is continuous and    $v'\in L^1(a, b)$  is finite almost everywhere, we must have
\begin{equation}
\label{psi<1}
  |\psi(v'(x))|< 1 \quad \text{for all }  x\in (a, b).
\end{equation}
This  implies that
$$
   v'_{|{(a, b)}} = \psi^{-1}(\psi(v')_{|{(a, b)}}) \in W^{1,1}_{{\rm loc}}(a,b)
$$
and it is decreasing. Therefore,  $v \in W^{2,1}_{{\rm loc}}(a, b)$ and it is concave.
\medskip

\noindent\texttt{Step 2}: ${u^j} =0$.
\\
Arguing by contradiction, we assume   that $u$ has a jump point at $z\in (a, b)$, and consider   the test function
$$
  \phi(x) :=
  \begin{cases}
   \frac{x-a}{z-a} &\hbox{if } x\in [a, z],\\[1ex]  0 & \hbox{if } x\in (z, b].
   \end{cases}
$$
Clearly, we have that
$$
D^s\phi= -\delta_z,
$$
where $\delta_z$  is the Dirac measure concentrated at $z$. Since $\phi(a)=\phi(b)=0$, and $|D^s\phi| =  \delta_z$ is absolutely continuous with respect to $|D^su|$ and its unique atom is $z$, it follows from \eqref{EE} that
\begin{align*}
\int_a^b \psi(D^au(x))  D^a\phi(x)\, dx - \int_a^b h(x)  \phi(x)\, dx   &= - \int_a^b\frac{D^su}{|D^su|}(x)  D^s\phi \\ & =   \int_a^b\frac{D^su}{|D^su|}(x) \, \delta_z
= \frac{D^su}{|D^su|} (z).
\end{align*}
On the other hand, integrating by parts,  it follows from \eqref{eqpsi} that
\begin{align}
\int_a^b \psi(D^au(x))  D^a\phi(x)\, dx &- \int_a^b h(x)  \phi(x)\, dx
 = \int_a^z \psi(D^au(x))  D^a\phi(x)\, dx - \int_a^z h(x)  \phi(x)\, dx
\\
&= \int_a^z \psi(v'(x))  \phi'(x)\, dx - \int_a^z h(x)  \phi(x)\, dx
\\
&=\psi(v'(z))  \phi(z) - \int_a^z \big(\psi(v'(x)) \big)' \phi(x)\, dx - \int_a^z h(x)  \phi(x)\, dx
\\
&=\psi(v'(z))  \phi(z) = \psi(v'(z)).
\end{align}
Hence, it follows that
\begin{align}
\label{v'z}
\psi(v'(z)) = \frac{D^su}{|D^su|} (z).
\end{align}
  Consequently, as the polar decomposition of measures (see, e.g., \cite[Cor. 1.29]{AmFuPa}) guarantees that
$$
 \left|\frac{D^su}{|D^su|}(z) \right| =1,
$$
we find from \eqref{v'z} that
$$
  |\psi(v'(z))| =1,
$$
which contradicts \eqref{psi<1}. Therefore, we conclude that $u^j =0$.
\medskip

\noindent\texttt{Step 3}: ${u^c} =0$.
\\
From the two previous steps, we already know that $u=u^a + u^c$ in $(a, b)$ and, hence, $u$ is continuous. Moreover, $u$ can be extended by continuity onto ${[a, b]}$. Let us prove that  $u$ is concave in ${[a, b]}$. On the contrary, assume that there exists an interval $[c, d] \subseteq{[a, b]}$ such that
$$
u(x) < u(c) + \tfrac{u(d) - u(c)}{d- c}(x-c) \quad \text{for all  }
x\in (c, d),
$$
 and consider the function $w\in BV(0,1)$ defined by
\begin{equation*}
w(x) =
\begin{cases}
\displaystyle u(c) + \tfrac{u(d) - u(c)}{d - c}(x-c) \  &\text{if } x\in [c, d],
\\[1ex] u(x) &\text{elsewhere.}
\end{cases}
\end{equation*}
It is clear that
$$
\int_0^1 \sqrt{1+|Dw|^2} < \int_0^1 \sqrt{1+|Du|^2}
$$
and, since  $w(x) > u(x) $ in $( c
, d
)$,
$$
 \int_0^1 h w\, dx \ge  \int_0^1 h u \, dx .
$$
Thus, we get
$$
 \mathcal{I}(w)=  \int_0^1 \sqrt{1+|Dw|^2} -  \int_0^1 h w \, dx <
\int_0^1 \sqrt{1+|Du|^2} -  \int_0^1 h u \, dx = \mathcal{I}(u),
$$
which contradicts the fact that $u$ is a global minimizer of the functional $\II$ defined by \eqref{M}.
Therefore,   $u$ being concave, it is locally Lipschitz  and, hence, ${u^c} =0$.

\smallskip
Since we have   proved that $u=u^a$, the remaining conclusions stated in (i) follow directly from the properties of $u^a$ established in Step 1. This ends the proof of the Assertion (i).
\medskip

\noindent \textit{Proof of Assertion \rm (iii).}
Thanks to (i) and (ii) we know that $u_{|(a,c)} \in W^{2,1}_{{\rm loc}}(a, c)\cap W^{1,1}(a, c)$ is concave   and $u_{|(c,b)} \in W^{2,1}_{{\rm loc}}(c,b)\cap W^{1,1}(c,b)$ is convex, and
$u^c=0$. Thus, as $u=v+u^j$, we have that, at the point $c$, either $u$ is  continuous, or it exhibits a jump discontinuity.
Arguing as in Step 1, we also see that $ \psi(v')\in  W^{1,1}(a, b)$,
\eqref{eqpsi} holds, and
\begin{equation}
  |\psi(v'(x))|< 1 \quad \text{for all }  x\in (a, b)\setminus \{c\}.
\end{equation}
If $|\psi(v'(c ))|<1$,   then
$$
|\psi(v')(x)| < 1 \quad \text{for all }  x\in (a,b),
$$
and hence the   argument used in Step 1 yields
 $v \in W^{2,1}_{{\rm loc}}(a,b)$.
Moreover, arguing as in Step 2 shows that $u^j=0$  and therefore  $u=v$. This proves that  $u\in W^{2,1}_{{\rm loc}}(a,b) \cap W^{1,1}(a,b)$.
\par
Suppose now  that   $|\psi(v'(c))| =1$.  Then, we have that  $\psi(v'(c)) =-1$, or equivalently $ v'(c) = -\infty$,  because $v_{|(a,c)}$ is concave and $v_{|(c,b)}$ is convex.
If $u$ is continuous at $c$,  it follows  that $u=v$ and, in particular,
$u'(c) =   -\infty$.
On the contrary, if $u$ has a jump discontinuity at $c$,   then  it follows   that
$u(c^-) >u(c^+) $, as the argument performed   in Step 2  of the proof of Assertion (i)   yields
\begin{align}
 \frac{D^su}{|D^su|} (c) =  \frac{D^ju}{|D^ju|} (c)  = \psi(v'(c)) = -1
\end{align}
and thus
$$
D^ju(c) = - |u(c^-)- u(c^+)| \, \delta_c.
$$
Finally, as  $u_{|(a,c)}= v_{|(a,c)}$ and $u_{|(c,b)}=v_{|(c,b)}+ u(c^+)- u(c^-)$,  we have that
$$
   u'(c^-) = u'(c^+) = v'(c) =-\infty.
$$
This concludes the proof.
\end{proof}

Essentially, Theorem \ref{th3.1} shows  that a bounded variation solution of equation \eqref{E} can  only  loose its regularity   at the endpoints, but not at the interior points, of  any interval where the right hand side of the equation  has a definite sign.
Our next two results, Theorem \ref{th3.2} and \ref{th3.3},  establish, on the contrary, the  complete regularity of a  bounded variation solution of  \eqref{E}. Precisely, Theorem \ref{th3.2}  guarantees the  regularity at the endpoints of any interval where  the sign of  the right hand side of the equation  is constant, by placing at   these points a suitable control on  its   decay rate to zero,  while Theorem \ref{th3.3} guarantees, instead, the regularity of a bounded variation solution at any interior point separating two adjacent   intervals  where the right hand side   changes sign, provided that a similar constraint on its decay rate to zero is imposed either on  the left, or on the right, of that point. Although this type of constraints were already introduced in some previous papers of ours, \cite{LGO-JDE20, LGO-ANS20, LGO-21}, here we will deliver a completely novel proof
of these regularity results based  on the strong maximum principle as expressed by Theorem \ref{th2.1}.
This  approach shows that the  condition  that we impose to earn regularity    is precisely  the assumption required by Theorem \ref{th2.1} so that the strong maximum principle holds for
  a certain equation satisfied by a shift of a local inverse of the considered solution.

\begin{theorem}
\label{th3.2}
Let $u$ be a bounded variation solution   of \eqref{E}.  Then, the following assertions  hold.
\begin{itemize}
\item[{\rm (j)}]
  If $f(x,u(x))\ge 0$   for almost all $x\in (a,b)$ and there exist  $\delta>0$ and $\mu \in L^1(a,a+\delta) $ such that
\begin{itemize}
\smallskip
\item[$\bullet$]  $ f(x,u(x))\le \mu(x) $   for almost all $x\in (a,a+\delta)$,
\item[$\bullet$]
$\displaystyle M(x) := \int_a^x\mu(t)\, dt>0$ for all $x\in (a, a+\delta]$,\,   and \,
$
\displaystyle
\int_{a}^{a+\delta} \frac{1}{\sqrt{M(x)}}\, dx = +\infty,
$
\end{itemize}
then $ u\in W^{2,1}_{\text{\rm loc}}[a,b)\cap W^{1,1}(a,b) $.
\smallskip
\item[{\rm (jj)}]
If $f(x,u(x))\ge 0$   for almost all $x\in (a,b)$ and there exist  $\delta>0$ and $\mu \in L^1(b-\delta,b) $ such that
\begin{itemize}
\smallskip
\item[$\bullet$]  $f(x,u(x))\le \mu(x)$ for almost all $x\in (b-\delta,b)$,
\item[$\bullet$] $\displaystyle M(x) := \int_x^b\mu(t)\, dt>0$ for all $x\in [b-\delta,b)$,\,  and \,
$
\displaystyle
\int_{b-\delta}^b \frac{1}{\sqrt{M(x)}}\, dx = +\infty,
$
\end{itemize}
then $ u\in W^{2,1}_{\text{\rm loc}}(a,b]\cap W^{1,1}(a,b)$.
\smallskip
\item[{\rm (jjj)}]
If $f(x,u(x))\le 0$   for almost all $x\in (a,b)$ and there    exist  $\delta>0$ and $\nu \in L^1(a,a+\delta) $ such that
\begin{itemize}
\smallskip
\item[$\bullet$]  $f(x,u(x))\ge \nu(x)$ for almost all $x\in (a,a+\delta)$,
\item[$\bullet$] $\displaystyle N(x) := \int_a^x\nu (t)\, dt<0$ for all $x\in (a, a+\delta]$,\,  and \,
$
\displaystyle
\int_{a}^{a+\delta} \frac{1}{\sqrt{-N(x)}}\, dx = +\infty,
$
\end{itemize}
then $ u\in W^{2,1}_{\text{\rm loc}}[a,b)\cap W^{1,1}(a,b) $.
\smallskip
\item[{\rm (jjjj)}]
If $f(x,u(x))\le 0$   for almost all $x\in (a,b)$ and there    exist  $\delta>0$ and $\nu \in L^1(b-\delta,b) $ such that
\begin{itemize}
\smallskip
\item[$\bullet$] $f(x,u(x))\ge \nu(x)$ for almost all $x\in (b-\delta,b)$,
\item[$\bullet$] $\displaystyle N(x) := \int_x^b\nu(t)\, dt<0$ for all $x\in [b-\delta,b) $,\,  and \,
$
\displaystyle
  \int_{b-\delta}^b \frac{1}{\sqrt{-N(x)}}\, dx = +\infty,
$
\end{itemize}
then $ u\in W^{2,1}_{\text{\rm loc}}(a,b]\cap W^{1,1}(a,b) $.
\end{itemize}
\end{theorem}

\begin{proof}
Let $u$ be a bounded variation solution of   equation \eqref{E}.
We will prove the validity of Assertion (j). As this proof can be easily adapted to complete
the proofs of the remaining assertions, these   are omitted here because repetitive.
By Theorem \ref{th3.1}(i), we already  know that   $u$ is concave, $ u\in W^{2,1}_{\text{\rm loc}}(a,b) \cap W^{1,1}(a,b) $, and it satisfies equation \eqref{E}   for almost all  $x\in (a, b)$.
Arguing by contradiction, we suppose that
\begin{equation}
\label{u'}
 u'(a^+)=+\infty.
 \end{equation}
 Then, there exists $c\in (a, a+\delta) $ such  that $u'(x) >0$ for all $x\in [a, c]$.   Thus, setting $\a= u(a )$ and  $\o= u(c)$,
we have that    $u: [a,c] \to [\a,\o]$ is invertible and the inverse function
$w: [\a,\o] \to[a,c]$
is continuously differentiable, with derivative
\begin{equation}
\label{w'}
  w'(t)=
\begin{cases}
\displaystyle
 \frac{1}{u'(w(t))} & \text{if } t\in (\a,\o],
 \\[1mm]
0 &  \text{if } t= \a.
\end{cases}
\end{equation}
By the concavity of $u$, $u'$ is  decreasing and hence $w'$ is  increasing. Thus,
\eqref{w'} implies that
\begin{equation}
\label{monw'}
  0<  w'(t) \leq w'(\o) = (u'(c))^{-2}\quad \hbox{for all}\;\; t \in (\a,\o].
\end{equation}
As $u' : (a,c] \to  [u'(c),+\infty)$ is locally absolutely continuous,  because
 $u\in W^{2,1}_{\text{\rm loc}}(a,b)$, and $w: (\a,\o] \to (a, c]$ is a  diffeomorphism of class $C^1$,  the composition
$$
  u'  \circ w : (\a,\o]   \to  [u'(c),+\infty)
$$
is locally absolutely continuous.  Hence,  the reciprocal function
$$
  \frac{1}{u'\circ w} : (\a,\o]   \to  \Big(0, \frac{1}{u'(c)}\Big]
$$
is locally absolutely continuous too.   Consequently, from \eqref{w'} we can infer  that $w' \in  W^{1,1}_{\text{\rm loc}}(\a,\o] $, with derivative
\begin{equation}
\label{w''}
w''(t)=\frac{-u''(w(t))}{(u'(w(t)))^2}\frac{1}{u'(w(t))} =-u''(w(t))(w'(t))^3 \quad
 \hbox{for almost all}\;\; t\in (\a,\o].
\end{equation}
Since the function $ u\in W^{2,1}_{\text{\rm loc}}(a,c]\cap W^{1,1}(a,c) $ solves the singular problem
\begin{equation}
\left\{ \begin{array}{l} -u''(x) = f(x,u(x)) \big(1+ (u'(x))^2\big)^\frac{3}{2}  \quad
\hbox{for almost all }  x\in (a, c),
\\[1mm]
u(a)=\a,\quad u'(a^+) =+\infty,
\end{array}
\right.
\end{equation}
and it satisfies the Lusin's $N$ property, it follows from \eqref{w'} and \eqref{w''} that
 the inverse function  $  w \in W^{2,1}_{\text{\rm loc}}(\a,\o] \cap C^{1}[\a,\o] $
is a solution of
\begin{equation}
\left\{ \begin{array}{l} w''(t)=f( w(t), t)  \big(  (w'(t))^2+1\big)^\frac{3}{2} \quad
\hbox{for almost all }  t\in (\a,\o),
\\[1mm] w(	\a)=a,\quad  w'(\a) =0.
\end{array}\right.
\end{equation}
 Therefore, if we define  $g: (\a, \o) \times \RR \times\RR \to \RR $ by
$$
g(t,s,\xi) := f(a + s, t) \big(1+\xi^2\big)^\frac{3}{2},
$$
and  $v\in W^{2,1}_{\text{\rm loc}}(\a,\o) \cap C^{1}[\a,\o] $ by
$$
   v:= w- a,
$$
the function $v$ is a non-trivial non-negative solution   of the initial value problem
\begin{equation*}
\left\{ \begin{array}{l} v''(t)= g(t,v(t), v'(t))
\quad \hbox{for almost all } t\in (\a, \o),
\\[1mm] v(\a)=0,\quad v'(\a) =0.
\end{array}\right.
\end{equation*}
Finally, let us introduce the function  $G:[0,\delta] \to \RR $   by  setting
$$
G(s) := M(a + s) \big(1+( u'(c))^{-2}\big)^\frac{3}{2},
$$
where $M$ has been defined in  (j).
From  the assumptions in (j) and the definitions of the functions  $g,  v,  w$ and $G$, we can infer from \eqref{monw'} that
 \begin{align}
 g(t, v(t), v'(t)) &= f(a+v(t), t) \big(1+( w'(t))^2  \big)^\frac{3}{2}  \\
 &  \le \mu(a+v(t)) \big(1+( u'(c))^{-2}\big)^\frac{3}{2} = G'(v(t)) \quad \text{for almost all } t\in (\a, \o),
\end{align}
and
$$
G(s) > 0 \; \text{ for all }  s \in (0,\d]  \qquad \text{and} \qquad   \int_{0}^{\delta} \frac{1}{\sqrt{G(s)}}\, ds = +\infty.
$$
Consequently, the conditions \eqref{0gG'} and \eqref{G>0} of (G) are both  fulfilled. Thus, thanks to Theorem \ref{th2.1}, it follows that $v'(\a)   >0$, which is a contradiction.
\end{proof}

The following  result establishes the regularity   of the solutions  at the interior points.

\begin{theorem}
\label{th3.3}
Let $u$ be a bounded variation solution of equation \eqref{E}.
 Then, the following
statements hold.
\begin{itemize}
\item[{\rm (h)}]
If there is $c\in (a,b)$ such that $f(x,u(x))\ge 0$   for almost all $x\in (a,c)$ and $f(x,u(x))\le 0$   for almost all $x\in (c,b)$ and either there  exist  $\delta>0$ and $\mu \in L^1(c-\delta,c) $ such that
\begin{itemize}
\item[$\bullet$] $f(x,u(x))\le \mu(x) $   for almost all $x\in (c-\delta,c)$,
\item[$\bullet$]
$ \displaystyle M(x) := \int_x^c\mu(t)\, dt>0$ for all $x\in [c-\delta,c) $,\,  and \,
$
\displaystyle
\int_{c-\delta}^c \frac{1}{\sqrt{M(x)}}\, dx = +\infty,
$
 \end{itemize}
 or there  exist  $\delta>0$ and $\nu \in L^1(c, c+\delta) $ such that
\begin{itemize}

\item[$\bullet$] $f(x,u(x))\ge \nu(x) $   for almost all $x\in (c, c+\delta)$,
\item[$\bullet$]
$ \displaystyle N(x) := \int_c^x\nu(t)\, dt<0$ for all $x\in (c, c+\delta] $,\, and  \,
$
\displaystyle
\int_c^{c+\delta}  \frac{1}{\sqrt{-N(x)}}\, dx = +\infty,
$
 \end{itemize}
then $ u\in W^{2,1}_{\text{\rm loc}}(a,b)\cap W^{1,1}(a,b) $.
\smallskip
\item[{\rm (hh)}]
If there is $c\in (a,b)$ such that $f(x,u(x))\le 0$   for almost all $x\in (a,c)$ and $f(x,u(x))\ge 0$   for almost all $x\in (c,b)$ and either there  exist  $\delta>0$ and $\nu \in L^1(c-\delta,c) $ such that
\begin{itemize}
\item[$\bullet$] $f(x,u(x))\ge \nu(x) $   for almost all $x\in (c-\delta,c)$,
\item[$\bullet$]
$ \displaystyle N(x) := \int_x^c\nu(t)\, dt<0$ for all $x\in [c-\delta,c) $,\, and \,
$
\displaystyle
\int_{c-\delta}^c \frac{1}{\sqrt{-N(x)}}\, dx = +\infty,
$
 \end{itemize}
 or there  exist  $\delta>0$ and $\mu \in L^1(c, c+\delta) $ such that
\begin{itemize}
\item[$\bullet$] $f(x,u(x))\le \mu(x) $   for almost all $x\in (c, c+\delta)$,
\item[$\bullet$]
$ \displaystyle M(x) := \int_c^x\nu(t)\, dt>0$ for all $x\in (c, c+\delta] $,\,  and \,
$
\displaystyle
\int_c^{c+\delta}  \frac{1}{\sqrt{M(x)}}\, dx = +\infty,
$
 \end{itemize}
then $ u\in W^{2,1}_{\text{\rm loc}}(a,b)\cap W^{1,1}(a,b) $.
\end{itemize}
\end{theorem}

\begin{proof}
Let $u$ be a bounded variation solution of  equation \eqref{E}, and, suppose, e.g.,
that the first alternative of the assumption (h) holds, i.e.,   there  exist  $\delta>0$ and $\mu \in L^1(c-\delta,c) $ such that
$$
f(x,u(x))\le \mu(x)     \;\; \hbox{for almost all}\;x\in (c-\delta,c), \qquad
M(x) := \int_x^c\mu(t)\, dt>0 \quad \hbox{for all}\; x\in [c-\delta,c),
$$
and
$$
\displaystyle
\int_{c-\delta}^c \frac{1}{\sqrt{M(x)}}\, dx = +\infty.
$$
By Theorem \ref{th3.1}, we already know that
either $ u\in W^{2,1}_{\text{\rm loc}}(a,b) \cap W^{1,1}(a,b)$,
or
$$
   u'(c^-)=-\infty= u'(c^+).
$$
On the other hand,  from  Theorem \ref{th3.2}(jj)  it follows that $ u\in W^{2,1}_{\text{\rm loc}}(a,c] $.  This rules out $u'(c^-)=-\infty $. Therefore, it follows that   $ u\in W^{2,1}_{\text{\rm loc}}(a,b) \cap W^{1,1}(a,b)$,  concluding the proof in this case. As the rest of the proof proceeds similarly, we omit the technical details.
\end{proof}

\section{Applications  to BVPs for the prescribed curvature equation}
\label{s4}

\noindent In this section we apply  the results   of Section \ref{s3} for establishing the regularity, up to the boundary, of the solutions of  non-autonomous one-dimensional  prescribed curvature equations subject to  possibly non-homogeneous  Dirichlet,  or Neumann,  or Robin boundary conditions, or else periodic boundary conditions. The case of   mixed   boundary conditions can be also dealt with, in an obvious way, by combining  these results and, hence, it is omitted here. In our exposition, we restrict ourselves to  discuss the regularity properties of a solution $u$ when $f(\cdot, u(\cdot))$ changes sign at most once. Other more general statements can be easily   inferred, in a  similar way, from Theorems \ref{th3.1}, \ref{th3.2} and \ref{th3.3}.

\subsection{The Dirichlet  problem}

Let us consider the    problem
\begin{equation}
 \label{DP}
\begin{cases}
\displaystyle
-\left( \frac{u'}{\sqrt{1+(u')^2}}\right)' = f(x,u), \quad 0<x<1,
\\
\displaystyle u(0) = \kappa_0, \; u(1) = \kappa_1 .
\end{cases}
\end{equation}
where  $f:(0,1)\times\RR \to \RR$
 and $\kappa_0, \kappa_1 \in\RR$ are  given.

\begin{dfn}
A function  $u\in BV(0,1)$  is a bounded variation solution of the Dirichlet problem  \eqref{DP} if
$f(\cdot , u(\cdot))\in L^1(0,1)$ and
\begin{align}
\int_0^1\frac{{D^au}(x) D^a\phi(x) }{\sqrt{1+(D^au(x))^2 }} \, dx+ \int_0^1 & \frac{D^su}{|D^su|}(x) \,   D^s\phi  + \sgn(u(0^+) - \kappa_0 )\phi(0^+)  \\
&   + \sgn(u(1^-) - \kappa_1) \phi(1^-)  = \int_0^1  f(x,u(x)) \phi(x) \, dx
\label{DPE}
\end{align}
for all $\phi \in BV(0,1) $ such that $|D^s\phi| $ is absolutely continuous with respect to $ |D^su|$,
$\phi(0^+)=0$ if $u(0^+)=\kappa_0$, and $\phi(1^-)=0$ if $u(1^-)=\kappa_1$.
\end{dfn}
We state below two sample regularity results for problem \eqref{DP}, assuming   that either
the right hand side of the equation in \eqref{DP}  does not change sign  in $(0,1)$, e.g., it is non-negative, or it   changes sign exactly once.

\begin{theorem}
\label{th4.1}
Let $u$ be a bounded variation solution of   \eqref{DP}. Assume that
\begin{itemize}
\item[(k)]    $f(x,u(x))\ge 0$   for almost all $x\in (0,1)$,
\item[(kk)] there  exist  $\delta>0$ and $\mu_0 \in L^1(0, \delta ) $ such that:
\begin{itemize}
\item[$\bullet$] $f(x,u(x))\le \mu_0(x) $   for almost all $x\in (0, \delta )$,
\item[$\bullet$]
$ \displaystyle M_0(x) := \int^x_{0} \mu_0(t)\, dt>0$ for all $x\in (0,\delta ] $,\, and  \,
$
\displaystyle
\int_{0}^{\delta } \frac{1}{\sqrt{M_0(x)}}\, dx = +\infty,
$
 \end{itemize}
\item[(kkk)] there  exist  $\delta >0$ and $\mu_1 \in L^1 \in (1-\delta ,1) $ such that
\begin{itemize}
\item[$\bullet$] $f(x,u(x))\le \mu_1(x) $   for almost all $x\in  (1-\delta ,1)$,
\item[$\bullet$]
$ \displaystyle M_1(x) := \int_x^{1} \mu_1(t)\, dt>0$ for all $x\in [1-\delta ,1) $,\, and \,
$
\displaystyle
\int_{1-\delta }^1 \frac{1}{\sqrt{M_1(x)}}\, dx = +\infty.
$
 \end{itemize}
  \end{itemize}
Then,  $ u\in W^{2,1} (0,1)$
and  it   satisfies both the   differential equation almost everywhere in $(0,1)$ and  the boundary conditions.
\end{theorem}

\begin{proof}
Theorems   \ref{th3.1} and \ref{th3.2}  guarantee that $ u\in W^{2,1}(0,1)$.
Testing  \eqref{DPE} against functions $\phi\in W^{1,1}(0,1)$ such that  $ \phi(0)=\phi(1)=0$ and integrating by parts, we infer that $u$ satisfies the    equation almost everywhere in $(0,1)$.   To show that $u$ fulfills the Dirichlet boundary conditions, we argue by contradiction.   Suppose, e.g., that $u(1) \neq\kappa_1$, and plug in \eqref{DPE}  a test function $\phi\in W^{1,1}(0,1)$ such that
$ \phi(0)=0$ and $ \phi(1)\neq  0$. Integrating by parts   in \eqref{DPE} yields
$$
\frac{u'(1) \,  \phi(1)  }{\sqrt{1+(u'(1))^2 }}
 + \sgn(u(1) - \kappa_1) \phi(1)  = 0
$$
and hence
$$
1 = |\sgn(u(1) - \kappa_1)| = \frac{|u'(1)|  }{\sqrt{1+(u'(1))^2 }},
$$
which is impossible, because $ u\in C^{1}[0,1]$. Therefore, $u(1) = \kappa_1 $. Similarly, one can prove that $u(0) = \kappa_0 $. This ends the proof.
\end{proof}

 Combining the proof of Theorem \ref{th4.1} with Theorem \ref{th3.3} yields the next result.

 \begin{theorem}
\label{th4.2}
Let $u$ be a bounded variation solution of  \eqref{DP}. Assume that
\begin{itemize}
\item[(l)]  there is $z\in (0,1)$ such that $f(x,u(x))\ge 0$   for almost all $x\in (0,z)$ and $f(x,u(x))\le 0$   for almost all $x\in (z,1)$,
\item[(ll)]
there  exist  $\delta>0$ and $\mu_0 \in L^1(0, \delta ) $ such that
\begin{itemize}
\item[$\bullet$] $f(x,u(x))\le \mu_0(x) $   for almost all $x\in (0, \delta )$,
\item[$\bullet$]
$ \displaystyle M_0(x):= \int^x_{0} \mu_0(t)\, dt>0$ for all $x\in (0,\delta ] $,\, and \,
$
\displaystyle
\int_{0}^{\delta } \frac{1}{\sqrt{M_0(x)}}\, dx = +\infty,
$
 \end{itemize}
\item[(lll)]
there  exist  $\delta >0$ and $\nu_1 \in L^1 \in (1-\delta ,1) $ such that
\begin{itemize}
\item[$\bullet$] $f(x,u(x))\ge \nu_1(x) $   for almost all $x\in  (1-\delta ,1)$,
\item[$\bullet$]
$ \displaystyle N_1(x) := \int_x^{1} \nu_1(t)\, dt<0$ for all $x\in [1-\delta ,1) $,\, and \,
$
\displaystyle
\int_{1-\delta }^1 \frac{1}{\sqrt{-N_1(x)}}\, dx = +\infty,
$
 \end{itemize}
\item[(llll)]
either there  exist  $\delta >0$ and $\mu_z \in L^1(z-\delta ,z) $ such that
\begin{itemize}
\item[$\bullet$] $f(x,u(x))\le \mu_z(x) $   for almost all $x\in (z-\delta ,z)$,
\item[$\bullet$]
$ \displaystyle M_z(x) := \int_x^z\mu(t)\, dt>0$ for all $x\in [z-\delta ,z) $,\, and \,
$
\displaystyle
\int_{z-\delta }^z \frac{1}{\sqrt{M_z(x)}}\, dx = +\infty,
$
 \end{itemize}
 or there  exist  $\delta >0$ and $\nu_z \in L^1(z, z+\delta ) $ such that
\begin{itemize}
\item[$\bullet$] $f(x,u(x))\ge \nu_z(x) $   for almost all $x\in (z, z+\delta )$,
\item[$\bullet$]
$ \displaystyle N_z(x) := \int_z^x\nu_z(t)\, dt<0$ for all $x\in (z, z+\delta ] $,\, and \,
$
\displaystyle
\int_z^{z+{\delta }}  \frac{1}{\sqrt{-N_z(x)}}\, dx = +\infty.
$
 \end{itemize}
  \end{itemize}
Then,  $ u\in W^{2,1} (0,1)$
and  it   satisfies both the   differential equation almost everywhere in $(0,1)$ and  the boundary conditions.
\end{theorem}

\begin{example}
 A simple situation where the assumptions of Theorems \ref{th4.1} and \ref{th4.2}
 are fulfilled is the  special case where
$$
   f(x,s) = h(x)\, k (s),
$$
with $h \in L^1(0,1)$ and $k\in C^0(\RR)$.   Suppose, for instance, that $u$  is a bounded variation solution of \eqref{DP} and that $k$ has a definite sign, e.g.,   $k(s)\ge0$ for all $s\in \mathrm{Range\,} u$,  and $h$ changes sign  once in $(0,1)$, e.g., there is $z\in (0,1)$ such that $h(x)> 0$   for almost all $x\in (0,z)$ and $h(x)<0$   for almost all $x\in (z,1)$.   Then, according to Theorem \ref{th4.2}, the regularity of $u$ is granted provided
$$
\int_{0}^{\frac{z}{2}} \left (\int_0^x  h(t) \, dt \right)^{-\frac{1}{2}}dx = +\infty, \qquad
\int_\frac{z+1}{2}^1 \left( \int_1^x h(t) \, dt \right)^{-\frac{1}{2}}dx = +\infty,
$$
and
$$ \text{either} \quad
\int^z_{\frac{z}{2}} \left (\int_x^z h(t) \, dt \right)^{-\frac{1}{2}}dx = +\infty, \quad
 \text{or}
 \quad
\int_z^\frac{z+1}{2} \left( \int_x^z h(t) \, dt \right)^{-\frac{1}{2}}dx = +\infty.
$$
 This and other similar regularity results, which can be easily deduced from Theorems \ref{th3.1}, \ref{th3.2} and  \ref{th3.3},  complete  the previous regularity results of  \cite{BHOO-JDE, BHOO-TS, OO-DIE, OO-JDE10, OO-DCDS, OS-NA} for \eqref{DP}.
\end{example}

\subsection{The  Neumann problem}

Let us consider the     problem
\begin{equation}
 \label{NP}
\begin{cases}
\displaystyle
-\left( \frac{u'}{\sqrt{1+(u')^2}}\right)' = f(x,u), \quad 0<x<1,
\\
\displaystyle\frac{u'(0)}{\sqrt{1+(u'(0))^2}} = \kappa_0, \quad  \frac{u'(1)}{\sqrt{1+(u'(1))^2}}= \kappa_1,
\end{cases}
\end{equation}
where
 $f:(0,1)\times\RR \to \RR$
 and $\kappa_0, \kappa_1 \in [-1,1]$ are given.

\begin{dfn}
A function  $u\in BV(0,1)$  is a bounded variation solution of   the Neumann problem  \eqref{NP} if
$f(\cdot , u(\cdot))\in L^1(0,1)$ and
\begin{equation}
\label{NPE}
\int_0^1\frac{u'(x) \phi'(x) }{\sqrt{1+(u'(x))^2 }} \, dx + \int_0^1  \frac{D^su}{|D^su|}(x) \,  D^s\phi
+\kappa_0 \ \phi(0^+) - \kappa_1  \phi(1^-)= \int_0^1  f(x,u(x)) \phi(x) \, dx
\end{equation}
for all $\phi \in BV(0,1) $ such that $|D^s\phi| $ is absolutely continuous with respect to $ |D^su|$.
\end{dfn}

 Here, we restrict ourselves  to    state    an illustrative    regularity result for  \eqref{NP}, under the assumption   that $f(\cdot, u(\cdot))$ changes sign   once in $(0,1)$.
Note that no assumptions are imposed on $f(\cdot, u(\cdot))$ at the boundary points, the regularity being guaranteed by  requiring that  $\kappa_0, \kappa_1 \in (-1,1)$,   which is optimal.

\begin{theorem}
\label{th4.3}
Let $u$ be a bounded variation solution   of \eqref{NP}.
Suppose  that conditions {\rm (l)} and {\rm(llll)} of  Theorem \ref{th4.2} hold and  $\kappa_0, \kappa_1 \in (-1,1)$.
Then, $ u\in W^{2,1} (0,1)$
and  it   satisfies both the   differential equation almost everywhere in $(0,1)$ and  the boundary conditions.
\end{theorem}

\begin{proof}
Theorems \ref{th3.1} and \ref{th3.3} guarantee that
$ u\in W^{2,1}_{\text{\rm loc}}(0,1)\cap W^{1,1}(0,1) $. We also  know  that
$$
\frac{u'}{\sqrt{1+(u')^2}} :  [0,1]  \to [-1,1]
$$
is continuous.   Plugging in \eqref{NPE} a function $\phi\in W^{1,1}(0,1)$ such that $ \phi(0)=\phi(1)=0$ and integrating by parts in $(-1,1)$, it becomes apparent that $u$ satisfies the    differential equation almost everywhere in $(0,1)$. Moreover, choosing  a test function $\phi\in W^{1,1}(0,1)$ such that   $ \phi(0)=0$ and $\phi(1)\neq 0$, and integrating by parts, we find that
$$
\frac{u'(1^-) \,  \phi(1)  }{\sqrt{1+(u'(1^-))^2 }}   - \kappa_1 \phi(1)  = 0
$$
and hence
$$
\frac{u'(1^-)  }{\sqrt{1+(u'(1^-))^2 }}  = \kappa_1.
$$
Similarly,   one can verify that
$$
\frac{u'(0^+)  }{\sqrt{1+(u'(0^+))^2 }}  = \kappa_0.
$$
Therefore,    $u$ satisfies the conormal boundary conditions.  Moreover, since $\kappa_0$, $\kappa_1 \in (-1,1)$,  we find that     $u \in C^1[0,1]$. Finally, the sign properties of $f(\cdot, u(\cdot))$ and, hence, of $u''$,  imply  that
\begin{equation}
\label{u''int}
\int_0^1  |u''(x)| \, dx = u'(0)- 2u'(z) +u'(1)\in \RR
\end{equation}
 and thus  $ u\in W^{2,1} (0,1)$. This ends the proof.
\end{proof}

Theorem \ref{th4.3} allows us to establish the regularity   of the solutions of \eqref{NP} whose existence was  proven in our preceding papers  \cite{LOR1, LOR2} (see also \cite{LGO-ANS19, LGO-JDE20, LGO-ANS20, LGO-21}).

\subsection{The Robin problem}

Let us consider the     problem
\begin{equation}
 \label{NR}
\begin{cases}
\displaystyle
-\left( \frac{u'}{\sqrt{1+(u')^2}}\right)' = f(x,u), \quad 0<x<1,
\\
 \displaystyle\frac{u'(0)}{\sqrt{1+(u'(0))^2}} + \la_0 u(0)  = \kappa_0, \quad  \frac{u'(1)}{\sqrt{1+(u'(1))^2}}
 + \la_1 u(1)= \kappa_1,
\end{cases}
\end{equation}
where  $f:(0,1)\times\RR \to \RR$
 and $\la_0, \la_1\in \RR\setminus\{0\}$, $ \kappa_0, \kappa_1 \in\RR$ are given.

\begin{dfn}
A function  $u\in BV(0,1)$  is a bounded variation solution of   the Robin problem  \eqref{NR} if
$f(\cdot , u(\cdot))\in L^1(0,1)$ and
\begin{align}
\label{NPR}
\int_0^1 \frac{u'(x) \phi'(x) }{\sqrt{1+(u'(x))^2 }} \, dx &+ \int_0^1  \frac{D^su}{|D^su|}(x) \,  D^s\phi
+  (\kappa_0  - \la_0 u(0^+) )\phi(0^+)  \\
&   +  (\la_1u(1^-) - \kappa_1) \phi(1^-)
= \int_0^1  f(x,u(x))  \phi(x) \, dx,
\end{align}
for all $\phi \in BV(0,1) $ such that $|D^s\phi| $ is absolutely continuous with respect to $ |D^su|$.
\end{dfn}
As above,  we just state here   an illustrative    regularity result for  problem \eqref{NPR} under the assumption   that the right hand side of the differential equation in \eqref{NPR}  changes sign   once in $(0,1)$.

\begin{theorem}
\label{th4.4}
Let $u$ be a bounded variation solution  of \eqref{NR}.
Suppose  that conditions {\rm (l)}, {\rm (ll)}, {\rm (lll)} and {\rm(llll)} of  Theorem \ref{th4.2} hold.
Then, $ u\in W^{2,1} (0,1)$
and  it   satisfies both the   differential equation almost everywhere in $(0,1)$ and  the boundary conditions.
\end{theorem}

\begin{proof}
Theorems \ref{th3.1}, \ref{th3.2}  and \ref{th3.3} guarantee that $ u\in W^{2,1} (0,1)$.  Moreover, plugging into \eqref{NPR} a test function $\phi\in W^{1,1}(0,1)$ such that $ \phi(0)=\phi(1)=0$,  and integrating by parts, it becomes apparent that $u$ satisfies the differential    equation almost everywhere in $(0,1)$,  while choosing  a test function $\phi\in W^{1,1}(0,1)$ such that  $ \phi(0)=0$ and $\phi(1)\neq 0$, and integrating by parts in $(0,1)$, it follows that
$$
 \frac{u'(1 ) \,  \phi(1)  }{\sqrt{1+(u'(1 ))^2 }} + (\la_1u(1 ) - \kappa_1) \phi(1)  = 0,
$$
and hence
$$
\frac{u'(1 )  }{\sqrt{1+(u'(1 ))^2 }}+  \la_1 u(1)     = \kappa_1.
$$
Similarly,  choosing  $\phi\in W^{1,1}(0,1)$ such that  $ \phi(0)\neq 0$ and $\phi(1)=0$ yields
$$
\frac{u'(0 )  }{\sqrt{1+(u'(0 ))^2 }} +\la_0u(0 ) = \kappa_0.
$$
 Therefore, $u$ satisfies the Robin boundary conditions. This concludes the proof.
\end{proof}

\subsection{The periodic   problem}

Let us consider the   problem
\begin{equation}
 \label{PP}
\begin{cases}
\displaystyle
-\left( \frac{u'}{\sqrt{1+(u')^2}}\right)' = f(x,u), \quad 0<x<1,
\\
u(0)=u(1),  \; u'(0) = u'(1),
\end{cases}
\end{equation}
where  $f:(0,1)\times\RR \to \RR$  is given.

\begin{dfn}
A function  $u\in BV(0,1)$  is a bounded variation solution of   problem \eqref{PP} if
$f(\cdot , u(\cdot))\in L^1(0,1)$ and
\begin{align}
\int_0^1\frac{D^a u(x) D^a\phi(x) }{\sqrt{1+(u'(x))^2 }} \, dx + \int_0^1  \frac{D^su}{|D^su|}(x) \,   D^s\phi &+\sgn(u(1^-)-u(0^+))(\phi(1^-) - \phi(0^+))
\\
\label{PPE}
&
= \int_0^1  f(x,u(x)) \phi(x) \, dx
\end{align}
for all $\phi \in BV(0,1) $ such that $|D^s\phi| $ is absolutely continuous with respect to $ |D^su|$  and $\phi(0^+)= \phi(1^-)$ if $u(0^+)=u(1^-) $.
\end{dfn}
Also in this case  we can state   a   regularity result for   problem \eqref{PP}  assuming  that  $f(\cdot,u(\cdot) )$ changes sign once in $(0,1)$.

\begin{theorem}
\label{th4.5}
Let $u$ be a bounded variation solution of  \eqref{PP}. Suppose that conditions {\rm (l)},  {\rm(llll)}, and either {\rm (ll)}, or  {\rm (lll)}, of Theorem \ref{th4.2} hold.
Then,  $ u\in W^{2,1} (0,1)$ and  it   satisfies both the   differential equation almost everywhere in $(0,1)$ and  the boundary conditions.
\end{theorem}

\begin{proof}
Assume that conditions {\rm (l)},  {\rm(llll)}  and, e.g.,   {\rm (lll)}  of Theorem \ref{th4.2} hold.
Then,  by Theorems   \ref{th3.1}, \ref{th3.2} and  \ref{th3.3},   $ u\in W^{2,1}_{\rm loc}(0,1]$ and,  in particular, $ u\in C^{1}_{\rm loc}(0,1]$.  Moreover, the function
$$
   \frac{u'}{\sqrt{1+(u')^2}} :  [0,1]  \to [-1,1]
$$
is continuous. Testing \eqref{PPE}  with functions  $\phi\in W^{1,1}(0,1)$ such that $ \phi(0) = \phi(1)=0$ and integrating by parts in $(0,1)$, it is easily seen that  $u$ satisfies the differential   equation almost everywhere in $(0,1)$.
Next, we choose  a test function $\phi\in W^{1,1}(0,1)$ such that $ \phi(0)=\phi(1)\neq 0$.  Then, integrating by parts yields
$$
\frac{u'(1) \,  \phi(1)  }{\sqrt{1+(u'(1))^2 }}  - \frac{u'(0^+) \,  \phi(0)  }{\sqrt{1+(u'(0^+))^2 }} =0,
$$
and hence
$$
\frac{u'(0^+) }{\sqrt{1+(u'(0^+))^2 }} =  \frac{u'(1)   }{\sqrt{1+(u'(1))^2 }}  \in (-1,1).
$$
This implies that
$
u'(0^+) = u'(1) \in \RR.
$
In particular, we infer that   $ u\in C^{1}[0,1]$. The sign properties of $f(\cdot, u(\cdot))$ and, hence, of $u''$  yields    \eqref{u''int}  and hence we can  conclude that  $ u\in W^{2,1} (0,1)$.
Finally,  arguing by contradiction,  suppose that $u(0)\neq u(1)$, and   pick  a test function $\phi\in W^{1,1}(0,1)$ such that $ \phi(1)-\phi(0)=1$. Integrating by parts, it follows  that
$$
\frac{u'(1) \,  \phi(1)  }{\sqrt{1+(u'(1))^2 }}  - \frac{u'(0) \,  \phi(0)  }{\sqrt{1+(u'(0))^2 }}
 + \sgn(u(1)-u(0))   = 0
$$
and hence, as  $u'(0) = u'(1)$,
$$
\frac{ |u'(1)|   }{\sqrt{1+(u'(1))^2 }}  =\frac{ |u'(1)|   }{\sqrt{1+(u'(1))^2 }}  |\phi(1)-\phi(0)|
 = |\sgn(u(1)-u(0)) |  = 1,
$$
which is impossible, because $ u\in C^{1}[0,1]$.   Therefore, we can conclude that  $u(0) = u(1) $. This ends the proof.
\end{proof}

 As far as concerns regularity, Theorem \ref{th4.5} complements and completes the existence results obtained   in  \cite{OO-JDE11, OO-CCM, OOR-NARWA, COZ-CCM, OO-OM}.

\end{document}